\documentclass[dvips,aap]{imsart}
\usepackage{amsthm,amsmath}
\RequirePackage[colorlinks,citecolor=blue,urlcolor=blue]{hyperref}

\startlocaldefs
\numberwithin{equation}{section}
\theoremstyle{plain}
\usepackage{amssymb}
\usepackage{pstricks,pst-fill,pst-node,pst-text,pst-plot,pst-coil}
\usepackage{graphicx}
\newcommand{\Leb}{\mathop{\rm Leb}\nolimits}
\newcommand{\W}{{\mathbf W}}

\newcommand{\cste}{\mathrm{Const.}}
\newcommand{\RR}{{\mathbb R}}
\newcommand{\xxx}{{\mathbf{X}}}

\newcommand{\Rup}{R^{\star}}
\newcommand{\e}{\varepsilon}
\newcommand{\Rdown}{R_{\star}}
\def\aire#1#2{\Leb_{#1}\left(#2\right)}
\def\vide{\varnothing}
\def\o{\mathrm{o}}
\def\O{\mathrm{O}}
\renewcommand{\P}[1]{\mathbb{P}\left(#1\right)}
\newcommand{\p}[1]{P\left(#1\right)}
\newcommand{\fexp}[1]{\exp\left(#1\right)}

\newenvironment{preuvesn}{\noindent {\sl Proof.~}}{\hfill\qed\par}
\newcommand{\leqb}{\Leb}
\newcommand{\Esp}{{\mathbb E}}

\def\d{\mathrm d}

\newrgbcolor{grey}{0.8 0.8 0.8}
\newtheorem{ansatz}{Ansatz}[section]
\newtheorem{lemma}{Lemma}[section]
\newtheorem{proposition}{Proposition}[section]
\newtheorem{remark}{Remark}[section]
\newtheorem{theorem}{Theorem}[section]
\newtheorem{definition}{Definition}[section]

\endlocaldefs

\begin{document}

\begin{frontmatter}

\title{Asymptotics of the visibility function in the Boolean model}
\runtitle{Visibility in the Boolean model}


\begin{aug}
\author{\fnms{Pierre} \snm{Calka}\corref{}\ead[label=e1]{pierre.calka@math-info.univ-paris5.fr}\thanksref{rem-pierre}}
\thankstext{rem-pierre}{Research partially supported by the French ANR project "mipomodim" No. ANR-05-BLAN-0017.}
\address{MAP5, U.F.R. de Math\'ematiques et Informatique\\
Universit\'e Paris Descartes\\
45, rue des Saints-P\`eres\\
75270  Paris Cedex 06, France.\\
\printead{e1}}
\affiliation{MAP 5, Universit\'e Paris Descartes}
\and
\author{\fnms{Julien} \snm{Michel}\ead[label=e2]{julien.michel@umpa.ens-lyon.fr}}
\address{Unit\'e de Math\'ematiques Pures\\
et Appliqu\'ees, UMR 5669\\
46 all\'ee d'Italie,\\
F-69364 Lyon Cedex 07, France.\\
\printead{e2}}
\affiliation{UMPA, \'Ecole Normale Sup\'erieure de Lyon}
\and
\author{\fnms{Sylvain} \snm{Porret-Blanc}\ead[label=e3]{sylvain.porret-blanc@umpa.ens-lyon.fr}}
\address{Unit\'e de Math\'ematiques Pures\\
et Appliqu\'ees, UMR 5669\\
46 all\'ee d'Italie,\\
F-69364 Lyon Cedex 07, France.\\
\printead{e3}}
\affiliation{UMPA, \'Ecole Normale Sup\'erieure de Lyon}

\runauthor{P. Calka et al.}
\end{aug}

\begin{abstract}
The aim of this paper is to give a precise estimate on the tail probability of the visibility
function in a germ-grain model: this function is defined as the length of the longest ray
starting at the origin that does not intersect an obstacle in a Boolean model.
We proceed
in two
or more dimensions using coverage techniques. Moreover, convergence results involving a
type I extreme
value distribution are shown in the two particular cases of small
obstacles or a large obstacle-free region.
\end{abstract}



\end{frontmatter}

\section{Presentation of the model and results}

In \cite{polya} G. P\'olya introduced the question of the visibility in a forest in a discrete lattice case as
well as in a random case. He first treated the problem of a person standing at the
origin of the regular square lattice of ${\mathbb R}^2$, when identical trees (discs with constant radius $R$) are situated
at the other points of the lattice. In this framework he showed that in order to see at a distance $r$
the radius $R$ should be (asymptotically when $r$ is large) taken as $1/r$. More recently V. Jankovi\'c
gave in \cite{jankovic} an elegant proof of a detailed version of this result. The random case
studied by G. P\'olya was the one of the visibility in one direction: we are here interested
in the global solution to this problem considering all directions
simultaneously. The spherical contact distribution which can be seen as the
infimum of the visibility over all directions has been intensively used for a geometric
description of random media (see e.g. \cite{muchestoyan, ratajsaxl,heinrich,lastschass,capassovilla,ballani}). 
In comparison, the total visibility,
{\it i.e.\/} the supremum of the same function,  
has been rarely
studied in the litterature. The work of reference is due to S. Zacks and is
strongly 
motivated by military applications (\cite{zacksLN},
see also his work with M. Yadin \cite{yadin}). 
However his interest was mainly focused on the probability that given points could be seen
and not on the total visibility. Very recently, the visibility problem has
been investigated in the hyperbolic disc by I. Benjamini et
alt. \cite{schramm}. In particular, the authors show the existence of a critical intensity
for the almost sure visibility at infinity. In this connection, one of the
consequences of our work will be that with probability one, we can see only at a
finite distance in the Euclidean space ${\mathbb R}^d$ (see Proposition \ref{highdim} and also Remark \ref{visibilityinfinity} concerning the
possibility to see at infinity).\\

In this paper, one of our goals is to present new distributional properties of
the total visibility in order to develop a future use of this indicator for
the study of porous media and more particularly in forestry. Potential
applications concern the optimization of directional logging of trees, the
measurement of competition level between growing trees in forest dynamics or
even an estimation of the light transmission through the canopy of a tree. In
such context, the total visibility seems to have an important role to play
even though it is understood that in some particular cases, another quantity of interest could be the mean of the visibility
in all directions.

The model is the following: consider a Boolean model (see \cite{molchanov,skm})
with random almost surely diameter-bounded {\em convex} grain $K$ with law $\mu$
based on a Poisson Point Process $\xxx$ with intensity measure the Lebesgue measure on ${\mathbb R}^d$, $d\geq 2$.
Define ${\mathfrak O}$ the {\em occupied phase} of this model,
$${\mathfrak O}=\bigcup_{x\in\xxx}(x\oplus K_{x}),$$
where $(K_{x})_{x\in\xxx}$ are independent identically distributed copies of $K$, independent of $\xxx$.
We condition this model by the event $O\notin {\mathfrak O}$ where $O$ is the
origin of ${\mathbb R}^d$. In particular, it has a positive probability equal
to $\exp(-\Esp[\leqb_{d}(K)])$.
We then define
the {\em visibility} in the following way:

\begin{definition}
Let ${\mathbf u}$ be a unit vector in $\RR^{d}$, the visibility in direction $\mathbf u$ is defined as
$$V({\mathbf u})=\inf\{r>0\ :\ r{\mathbf u}\in{\mathfrak O}\},$$
the total visibility is defined as
$${\mathfrak V}=\sup_{\|{\mathbf u}\|=1}V({\mathbf u}).$$
\end{definition}

\vspace*{18mm}
\begin{center}
\psset{unit=8.5mm}
\pscircle[linecolor=black,fillcolor=black,fillstyle=solid](0,0.2){0.04}
\rput[b](-0.3,0.1){$O$}
\pscircle[linecolor=grey,fillcolor=grey,fillstyle=solid](1.21,2.113){0.32}
\pscircle[linecolor=grey,fillcolor=grey,fillstyle=solid](0.123,1.345){0.22}
\pscircle[linecolor=grey,fillcolor=grey,fillstyle=solid](-1.564,0.234){0.12}
\pscircle[linecolor=grey,fillcolor=grey,fillstyle=solid](-0.534,-0.765){0.54}
\pscircle[linecolor=grey,fillcolor=grey,fillstyle=solid](0.52,1.53){0.09}
\pscircle[linecolor=grey,fillcolor=grey,fillstyle=solid](0.23,-1.24){0.24}
\pscircle[linecolor=grey,fillcolor=grey,fillstyle=solid](-0.978,-0.142){0.16}
\pscircle[linecolor=grey,fillcolor=grey,fillstyle=solid](-1.12,1.45){0.36}
\pscircle[linecolor=grey,fillcolor=grey,fillstyle=solid](1.987,0.21){0.12}
\pscircle[linecolor=grey,fillcolor=grey,fillstyle=solid](1.03,-1.13){0.46}
\pscircle[linecolor=grey,fillcolor=grey,fillstyle=solid](0.0214,-1.2){0.15}
\psline[linecolor=black,linewidth=1pt]{->}(0,0.2)(0.3,-0.1)
\rput[b](0.3,0.1){\small${\mathbf u}$}
\psline[linecolor=black,linestyle=dashed,linewidth=0.5pt]{-*}(0,0.2)(0.88,-0.68)
\rput[b](1.2,-0.65){\small$V({\mathbf u})$}
\end{center}
\vspace*{11mm}
\begin{figure}[h]
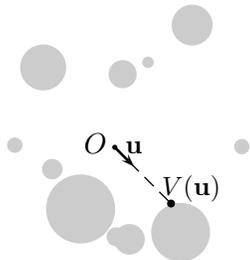

\caption{Directional visibility}\label{dessin-foret}
\end{figure}

For a convex body $K$ and ${\mathbf u}\in {\mathbb S}^{d-1}$, we define the width of $K$
in the direction $v$ orthogonal to ${\mathbf u}$ as
$$W_u({K})=\sup_{x,y\in {K}} \langle y-x, {\bf v}\rangle.$$
The mean width of $K$ is denoted by $\W(K)$ and satisfies
$$\W(K)=\frac{1}{d\omega_d}\int W_u(K)\,\d\sigma_d(u),$$
where $\sigma_d$ is the uniform measure on ${\mathbb S}^{d-1}$ and $\omega_d$
is the $d$-dimensional Lebesgue measure of the unit-ball in ${\mathbb R}^d$.

It is well known that the law of the directional visibility is exponential, indeed
\begin{lemma}\label{une-direction}
For each unit direction ${\mathbf u}$ in $\RR^{d}$ one has
$$\P{V({\mathbf u})>r}=\exp(- r\Esp[\W(K)]).$$
\end{lemma}

The aim of this paper is to give precise estimates on the tail probability for the visibility
function in all directions: in section 2 we present the general method of coverage processes
used throughout this paper. In section 3 we give sharp upper and lower bounds in the same exponential
order $\exp(-\cste r)$ in dimension two and in the two cases of circular obstacles
and of more general rotation-invariant random obstacles.
In 
higher dimensions results on coverage processes are more sparse, thus the bounds presented in section 4 for 
dimension $d\geq 3$ are rougher.
Section 5 is devoted to two similar convergence results
for the asymptotics of the visibility with small obstacles, and when the spherical contact length is conditionned to be large. Both results
state a convergence in law towards a Gumbel distribution, are valid for any
dimension $d\geq 2$ and are based on an extension of a result of
\cite{janson}. 

The present work has been first announced in a small note \cite{CMPB09}.
\section{Random coverage of the circle and the sphere and visibility}

The visibility up to length $r$ may be blocked only by those obstacles that intersect $B_{d}(0,r)$, the number $N_{r}$ of
those obstacles is Poisson distributed with parameter 
$$g_{d}(r,K)=\Esp[\Leb_{d}((B_{d}(0,r)\oplus \check{K})\setminus \check{K})],$$
where $\check{K}=\{-x\,:\,x\in K\}$,
and those obstacles $K_{i}$, $i\in\{1,\ldots,N_{r}\}$, are independent and identically 
distributed with the same law. Each of them projects a {\em shadow} of solid angle $S_{i}$
on the sphere $r{\mathbb S}^{d-1}=\partial B_{d}(0,r)$, which is a random cap $C_{i}$: 
\begin{equation}
C_{i}=\{(r,{\mathbf u})\in{\mathbb R}_{+}\times{\mathbb S}^{d-1}\ :\ \exists s\leq r
\hbox{\ s.t.\ } (r,{\mathbf u})\in K_{i}\},\label{cap}
\end{equation}
and the solid angle $S_{i}$ is equal to
\begin{equation}
S_{i}=\{{\mathbf u}\in{\mathbb S}^{d-1}\ :\ (r,{\mathbf u})\in C_{i}\}.\label{angle-solide}
\end{equation}
We have thus the following ansatz:
\begin{ansatz}\label{coverage}
The visibility ${\mathfrak V}$ is greater than $r$ if and only if the solid angles $S_{i}$, $i\in\{1,\ldots,N_{r}\}$
do not cover the sphere  ${\mathbb S}^{d-1}$.
\end{ansatz}
This equivalence links coverage properties of the sphere with our initial problem: this problem of random coverings 
has been quite intensively studied in the literature in the two-dimensional case, see for instance \cite{dvoretzky,kahane,fan} 
in the context of Dvoretzky covering, and 
\cite{siegelholst,shepp,stevens,yadin} for a more general approach. The properties of coverings in
higher dimensions are less known, let us cite the works of S. Janson (\cite{janson} and other papers) dealing with
some asymptotic properties for coverage processes with small caps.\\

Let us remark that each obstacle $K_i$ is distributed according to the law
$$\displaystyle g_d(r,K)^{-1}{\bf 1}_{(x\oplus K)\cap B_d(0,r)\ne \emptyset}\,\d x
\d\mu(K).$$ We denote by $\tilde{\nu}_r$ the distribution of the associated
solid angle $S_i$. The probability measure $\tilde{\nu}_r$ is naturally
invariant under rotations. We call $P(\tilde{\nu}_r,n)$ the probability that
$n$ independent solid angles distributed as $\tilde{\nu}_r$ cover the sphere. 
With Ansatz \ref{coverage} the following result becomes straightforward:
\begin{proposition}\label{equiv-recouvre}
For every $r>0$, 
\begin{equation}
  \label{lawwithcover}
\P{{\mathfrak V}\geq r}=\fexp{-g_{d}(r,K)}\sum_{n\geq 0}\frac{(g_{d}(r,K))^n}{n!}(1-\p{\tilde{\nu}_r,n}).
\end{equation}
\end{proposition}
When in dimension two, for every probability measure $\nu$ on $[0,1]$,
$P(\nu, n)$ denotes the covering probability of the circle with perimeter one
by $n$ i.i.d. isotropic arcs of $\nu$-distributed length. 
We shall denote by $\nu_{r}$ the probability
law of the random arc shadowed by an obstacle intersecting the disc $B_{2}(0,r)$.\\
When $K=B_{2}(0,R)$ for a constant radius $R$, one has for instance the following values for the characteristic
data of proposition \ref{equiv-recouvre}:
\begin{eqnarray}
\label{lengthmeasure-1}g_{2}(r,K)&=&\pi(2rR+r^2),\\
\frac{\d\nu_r}{\d u}(u)&=&\frac{\pi r}{rR+\frac{r^2}{2}}\,{\bf
  1}_{\left[0,\frac{1}{\pi}\arctan\left(\frac{R}{r}\right)\right]}(u)\times\nonumber\\
  &&\left(r\sin(2\pi
    u)+\frac{\sin(\pi u)(R^2+r^2\cos(2\pi u))}{\sqrt{R^2-r^2\sin^2(\pi
        u)}}\right)\nonumber\\
\label{lengthmeasure-2}&&+\frac{\pi R^2}{rR+\frac{r^2}{2}}\,{\bf
    1}_{\left[\frac{1}{\pi}\arctan\left(\frac{R}{r}\right),\frac{1}{2}\right]}(u)
    \,\frac{\cos(\pi
    u)}{\sin^3(\pi u)}.
\end{eqnarray}
\begin{remark}[The two-dimensional case]
Let us recall that the covering probability $\p{\nu,n}$ has been computed in dimesion two by
Stevens in the case of arcs with deterministic lengths \cite{stevens} and Siegel \& Holst in
the general case \cite{siegelholst}: for every probability measure $\nu$ on
$[0,1]$ and $n\in {\mathbb N}^*$, 
\begin{equation}
  \label{sh}
\p{\nu,n}
=\sum_{k=0}^n (-1)^{k}\binom{n}{k}
\int \left[\prod_{i=1}^k F_{\nu}(u_i)\right]\left[\sum_{i=1}^k \int_0^{u_i}
  F_{\nu}(t)\,\d t\right]^{n-k} \,\d{\tilde{\lambda}}_k(u)
\end{equation}
where $F_{\nu}$ is the cumulative distribution function of $\nu$ and $\tilde{\lambda}_k$,
$k\in {\mathbb N}^*$,  is
the normalised uniform measure on the simplex $\{(x_1,\cdots,x_k)\in
[0,1]^k;x_1+\cdots+x_k=1\}$.

We can consequently substitute $\p{\nu,n}$
by its expression (\ref{sh}) to get an explicit formula for $\P{{\mathfrak V}\geq r}$,
$r>0$. Nevertheless, it seems more or less intractable for doing some asymptotic
estimations. 
\end{remark}

\section{Sharp asymptotics in dimension two}

In dimension two it becomes possible to give sharp estimates for the tail probability, for instance we have
from lemma \ref{une-direction} the following lower bound:
\begin{equation}
\P{{\mathfrak V}>r}\geq \P{V_{\mathbf u}>r}=\exp(-r\Esp[\W(K)]),\label{mino-triviale}
\end{equation}
in this section we shall give two sharper lower bounds and an upper bound:
\begin{itemize}
\item When obstacles are fixed discs with constant radius $R$, 
for each $\mu\in(0,2/R)$ there exists a function $\varepsilon_{\mu}$ converging to $0$ as $r\to\infty$ 
such that:
\begin{equation}
\P{\mathfrak{V}\geq r} \geq \mu r\exp({-2Rr})\left(1+\varepsilon_{\mu}(r)\right),\label{mino-main}
\end{equation}
for $r$ large enough;
\item When obstacles are discs with random radius $R$, we have for $r$ large enough
\begin{eqnarray}
\P{{\mathfrak V}\geq r}&\leq& 2(\pi(2r\Esp[R]+r^2)+2)\times\nonumber\\
&&\fexp{-\pi(2r\Esp[R]+r^2)m_r},\label{eq:10}\\
\P{{\mathfrak V}\geq r}&\geq&  2\pi m_r (2r\Esp[R]+r^2)\fexp{-\pi m_r(2r\Esp[R]+r^2)}+\nonumber\\
&&\fexp{-2\pi m_r(2r\Esp[R]+r^2)},\label{eq:10.1}
\end{eqnarray}
where $m_{r}$ is the mean of law $\nu_{r}$, satisfying
\begin{equation}m_{r}\sim \frac{2\Esp[R]}{\pi r}.\label{equiv-moyenne}
\end{equation}
>From (\ref{eq:10}), (\ref{eq:10.1}) and (\ref{equiv-moyenne}) we obtain directly the following theorem:
\begin{theorem}\label{logestimate-R}
For the Boolean model in dimension two with random discs,
the asymptotics of the visibility is given by
$$\lim_{r\to+\infty}\frac{1}{r}\,\log \P{{\mathfrak V}\geq r}=-2\Esp[R].$$
\end{theorem}
\end{itemize}
The estimate (\ref{mino-main}) is obtained
without the use of the previous section, by more
direct considerations based on the estimation of the visibility in a finite
number of directions. Details of the proof are postponed to the appendix 1.

In the following lines, we shall prove estimates (\ref{eq:10}) and (\ref{eq:10.1}), as well as 
a generalisation of theorem \ref{logestimate-R} for
generall convex shapes. Most of the arguments rely on comparison results for covering
probabilities of the circle. 

\subsection{Lower bound via comparison of coverage probabilities, estimate (\ref{eq:10})}

In this subsection and the next we shall follow the main steps of a
previous work \cite{calka03} that dealt with the circumscribed radius of the typical 
Poisson-Voronoi cell: a comparison result on covering probabilities states
that if two probability measures (the laws of the arcs) are comparable in some sense, then the
coverage probabilities are also comparable. Up to now the ordering induced by
the concentration around
the mean has been the main criterion for comparing covering probabilities
\cite{shepp, calka03} but the convex ordering (which is implied by the previous ordering)
is in fact enough to deduce the required inequalities.

For the special case of random discs, the computation of the parameters of equation (\ref{lawwithcover}) is done
in the following way. Let us denote by $\mu$ the law of the radius of the random discs of the Boolean model.
The law of the normalised lengths of the shadowed arcs on the circle is denoted by $\nu_{r}$, 
this law is the image of the couple $(R,U)$ with law $\mu\otimes{\mathcal U}(0,1)$
by the map $(R,U)\mapsto \ell$:
\begin{eqnarray}
\ell=\left\{\begin{array}{ll}
\displaystyle{\frac{1}{\pi}\arcsin\frac{1}{\sqrt{1+\frac{r}{R}\left(2+\frac{r}{R}\right)\,U}}}&
    \mbox{ if $0\leq U\leq \dfrac{r}{r+2R}$;}\\
    \displaystyle{\frac{1}{\pi}\arccos\frac{\frac{r}{R}+\left(2+\frac{r}{R}
    \right)\,U}{2\sqrt{1+\frac{r}{R}\left(2+\frac{r}{R}\right)\,U}}}&\mbox{
    if $\dfrac{r}{r+2R}\leq U\leq 1.$}\end{array}\right.  \label{length2}
\end{eqnarray}
As $r$ tends to infinity, we obtain 
$\ell\simeq({1}/{\pi})\arcsin({R}/{(r\sqrt{U})})$ so that the asymptotics of the expectation becomes
$m_{r}:=\Esp[\ell]\simeq{2\Esp[R]}/({\pi r})$.\\

We shall use here the convex domination
of measures, let us first remark the following:
for any probability measure $\nu$ on $[0,1/2]$ with mean $m$ and any convex
function $f:[0,1/2]\longrightarrow {\mathbb R}$, we have
\begin{equation}
  \label{convexorder}
\int f\,\d\nu\leq (1-2m)f(0)+ 2m f(1/2),
\end{equation}
which means that $\nu<_{cv} \left[(1-2m)\delta_0+(2m)\delta_{1/2}\right]$
(where $<_{cv}$ denotes the usual convex order \cite{mullerstoyan}).\\

It is a consequence of the proof of theorem 13 in \cite{calka03} that the convex order
implies the order of the covering probabilities: if $\mu_1$ and $\mu_2$ are two
probability measures on $[0,1/2]$ such that $\mu_1<_{cv}\mu_2$ then
$\p{\mu_1,n}\leq \p{\mu_2,n}$. 
Inserting the inequality $P(\nu_r,n)\le P((1-2m_r)+2m_r\delta_{1/2},n)$  in
(\ref{lawwithcover}) for every $n\in{\mathbb N}$, we get
that
\begin{eqnarray}
\P{{\mathfrak V}\geq r}&\geq& \fexp{-\pi(2r\Esp[R]+r^2)}\times\nonumber\\
&&\sum_{n\geq
  0}\!\frac{(\pi(2r\Esp[R]+r^2))^n}{n!}(1-\p{(1-2m_r)\delta_{0}+2m_r\delta_{1/2},n}).\label{eq:123}
\end{eqnarray}
It follows from (\cite{calka03}, Corollary 1) that 
$$1-\p{(1-2m_r)\delta_{0}+2m_r\delta_{1/2},n}=2nm_r(1-m_r)^{n-1}+(1-2m_r)^{n}.$$
Inserting that result in (\ref{eq:123}), we obtain that
\begin{eqnarray}
\P{{\mathfrak V}\geq r}&\geq&  2\pi m_r (2r\Esp{[R]}+r^2)\fexp{-\pi m_r(2r\Esp[R]+r^2)}+\nonumber\\
&&\fexp{-2\pi m_r(2r\Esp[R]+r^2)}.
\label{liminf}
\end{eqnarray}

\subsection{Upper bound via comparison of coverage probabilities, estimate (\ref{eq:10.1})}

By Jensen's inequality, we have that $\nu_r>_{cv}\delta_{m_r}$. Consequently,
the same argument as for the lower bound shows that
$P(\nu_r,n)\geq P(\delta_{m_r},n)$ for every $n\in{\mathbb N}$.
Inserting this inequality in (\ref{lawwithcover}), we have
\begin{equation}\label{majorationmoyenne}
\P{{\mathfrak V}\geq r}\leq \fexp{-\pi(2r\Esp[R]+r^2)}\sum_{n\geq
  0}\frac{(\pi(2r\Esp[R]+r^2))^n}{n!}(1-\p{\delta_{m_r},n}).
\end{equation}

For the estimation of $(1-\p{\delta_a,n)}$, Shepp obtained a basic inequality \cite{shepp}
which holds for $a\in[0,1/4]$ and
$n\in{\mathbb N}^*$ and is easier to use than Steven's explicit formula:
\begin{equation}
  \label{eq:shepp}
1-\p{\delta_a,n}\leq \frac{2(1-a)^{2n}}{\int_0^{a} (1-a-t)^{n}\,\d t +(\frac{1}{4}-
a)(1-2a)^{n}}.
\end{equation}
A straightforward consequence of (\ref{eq:shepp}) is that for every $a\in
[0,1/4]$ and $n\in {\mathbb N}$, we have 
\begin{equation}
  \label{eq:sheppameliore}
1-\p{\delta_{a},n}\leq 2(n+1)(1-a)^{n-1}.
\end{equation}
In particular, for $r$ sufficiently large, the mean $m_r$ is in the interval $[0,1/4]$
so the equality (\ref{lawwithcover}) combined with (\ref{eq:sheppameliore})
leads us to 
$$\P{{\mathfrak V}\geq r}\leq 2(\pi(2r\Esp[R]+r^2)+2)\exp(-\pi(2r\Esp[R]+r^2)m_r).$$
It remains to use the estimation on the mean $m_{r}$ to get that
\begin{equation}
  \label{limsup}
  \limsup_{r\to+\infty}\frac{1}{r}\,\log \P{{\mathfrak V}\geq r}\leq -2\Esp[R].  
\end{equation}
(\ref{liminf}) and (\ref{limsup}) now complete the proof of Theorem \ref{logestimate-R} for random
discs.

\subsection{The case of general convex shapes}

Let ${K}$ be a random convex body of ${\mathbb R}^2$ containing the origin, which is supposed to be invariant
under any rotation and is such that its diameter is bounded almost surely by a
constant $D>0$. For instance, ${\mathcal K}$ can be the image of a
deterministic convex body by a uniform random rotation. 
%
By the rotation-invariance of ${K}$, we have $\Esp[W_u({
  K})]=\Esp[{\bf W}({K})]$. 
\begin{theorem}\label{logestimate}
For the Boolean model with random rotation-invariant grains distributed as
${K}$, the asymptotics of the visibility is given by
$$\lim_{r\to+\infty}\frac{1}{r}\,\log \P{{\mathfrak V}\geq r}=-\Esp[{\bf W}({
  K})].$$
  \end{theorem}
  \begin{preuvesn}
The lower bound is obtained by lemma \ref{une-direction}. It remains to show that 
$$\limsup_{r\to +\infty}\frac{1}{r}\log \P{{\mathfrak V}\geq r}\leq -\Esp[{\bf W}({
  K})].$$
In order to do it, we need an intermediary geometric lemma whose proof is
postponed to the appendix 2:
\begin{lemma}\label{geometriclemma}
Let  $L$ be a convex body of diameter bounded by $D$ containing the origin. We define
$\Psi(r{\bf u}+L)$ as the angle of vision of $(r{\bf u}+L)$ from
$O$. Then
$$\lim_{r\to +\infty}r\Psi(r{\bf u}+L)=W_u(L),$$
the limit being uniform over all unit-vectors ${\bf u}$ and all such convex bodies $L$. 
\end{lemma}
As in the proof of Proposition \ref{equiv-recouvre}, the event $\{{\mathfrak
  V}\geq r\}$, $r>0$, can be seen as the uncovering of the circle $C(0,r)=r{\mathbb S}^1$ by
the 'shadows' produced by the obstacles $(x\oplus{ K}_x)$ such that $(x\oplus{
  K}_x)\cap B_2(O,r)\ne \emptyset$. 

Let $\varepsilon>0$. By Lemma
\ref{geometriclemma}, let us fix $r_{\varepsilon}>0$ such that for every $x$
such that $\|x\|>r_{\varepsilon}$ and every convex body $L$ (with a diameter
bounded by $D$), we have 
\begin{equation}
  \label{eq:epsilonpres}
\|x\|\Psi(x+L)\geq (W_{x/\|x\|}(L)-\varepsilon).  
\end{equation}
Then for $r>r_{\varepsilon}+D$, the probability of uncovering the circle
$C(0,r)$ is greater if we only keep the shadows produced by the
obstacles $(x\oplus{K}_x)$ such that $r_{\varepsilon}<\|x\|<r-D$. In that
case, such a shadow is a random
rotation-invariant arc on the
circle $C(O,r)$ whose normalised length is $(2\pi)^{-1}\Psi(x\oplus{K}_x)$. Let us denote
by $\eta_r$ the mean of $(2\pi)^{-1}\Psi(Z\oplus{K}_Z)$ when $Z$ is uniformly distributed
in $B_2(r-D)\setminus B_2(r_{\varepsilon})$ and ${K}_Z$ is independent from $Z$
and distributed as ${K}$. Following the method already used to obtain the
upper-bound (\ref{majorationmoyenne}), we have
$$
\P{{\mathfrak V}\geq r}\leq
e^{-\pi((r-D)^2-r_{\varepsilon}^2)}\sum_{n=0}^{+\infty}\frac{\pi((r-D)^2-r_{\varepsilon}^2)}{n!}
(1-P(\delta_{\eta_r},n)).$$
Forecasting that $\eta_r$ will be small enough, we may apply the inequality
(\ref{eq:sheppameliore}) in order to obtain that
\begin{equation}
  \label{eq:102}
\P{{\mathfrak V}\geq r}\leq 2(\pi((r-D)^2-r_{\varepsilon}^2)+2)\exp(-\pi((r-D)^2-r_{\varepsilon}^2)\eta_r).
\end{equation}
Let us now estimate the mean $\eta_r$: using (\ref{eq:epsilonpres}), we get
\begin{eqnarray*}
\eta_r&=&\frac{1}{2\pi}E_{{K}}\left[\int_0^{2\pi}\int_{r_{\varepsilon}}^{r-D}\Psi(\rho
  u_{\theta}\oplus{K})\rho \frac{\d\rho}{\pi((r-D)^2-r_{\varepsilon}^2)}\,
  \d\theta\right]\\
&\geq& \frac{1}{2\pi}E_{{K}}\left[\int_0^{2\pi}(W_{\theta}({
    K})-\varepsilon)\,\d\theta\right]\cdot\int_{r_{\varepsilon}}^{r-D}\frac{\d\rho}{\pi((r-D)^2-r_{\varepsilon}^2)}\\
&= & (\Esp[{\bf
W}({
K})]-\varepsilon)\int_{r_{\varepsilon}}^{r-D}\frac{\d\rho}{\pi((r-D)^2-r_{\varepsilon}^2)}\\
&\underset{r\to +\infty}{\sim} &\frac{\Esp[{\bf
W}({
K})]-\varepsilon}{\pi r}.
\end{eqnarray*}
Inserting this last result in (\ref{eq:102}), we obtain that
$$\limsup_{r\to +\infty}\frac{1}{r}\log \P{{\mathfrak V}\geq r}\leq -(\Esp[{\bf
  W}({K})]-\varepsilon).$$
When $\varepsilon$ goes to $0$, we get the required result.
  \end{preuvesn}

\section{Rough estimates in dimension greater than three}

The problem of maximal visibility in a Boolean model is investigated in
${\mathbb R}^d$ with deterministic radii $R_x=R$, $x\in \xxx$. The obstacles are balls of deterministic radius $R$.
The same connection between the distribution of ${\mathfrak V}$ and the non-covering of the sphere by random circular 
caps occurs. What prevents us from obtaining the analogue of Theorem \ref{logestimate} is that the calculation of the 
probability to cover the sphere with caps of random radii is not known. We have to restrict ourselves to coverings of the 
sphere with caps of deterministic radii. This explains that the following result is weaker than Theorem \ref{logestimate}:
\begin{proposition}\label{highdim}
In dimension $d\geq 3$, we have
$$\liminf_{r\to +\infty}\frac{1}{r}\,\log\P{{\mathfrak V} \geq r}\geq -\omega_{d-1}R^{d-1}$$
and
$$ \limsup_{r\to +\infty}\frac{1}{r}\,\log \P{{\mathfrak V}\geq r}\leq
  -\frac{1}{d}\omega_{d-1} R^{d-1}.$$
\end{proposition}
\begin{preuvesn}
As in the two-dimensional case, the lower-bound is obtained by considering the visibility in a fixed direction.

Let us focus on the upper-bound: the maximal visibility is larger than $r>0$ if and only if the shadows produced 
by the obstacles on the sphere centred at the origin and of radius $r$ do not cover that sphere. The concerned balls are 
those such that their centres are at distance $\rho\in [R,R+r]$ from the origin. Since we look for an upper-bound of a 
probability of non-covering of the sphere by random circular caps, we can take less and smaller caps. For sake of simplicity, we only keep
 the shadows produced by the balls with a centre at distance $\rho\in [R, r]$.

For such a ball, it comes from (\ref{length2}) that the angular radius of its shadow on the sphere is at least 
$\arcsin\left(\frac{R}{r}\right)$, which is bigger than $\frac{R}{r}$.

In conclusion, the probability that the maximal visibility is greater than $r$
is lesser than the probability of non-covering of the unit-sphere by  a
Poissonian number (of mean $\omega_d(r^d-R^d)$) of circular caps of angular
radius $\frac{R}{r} $. Upper-bounds for covering probabilities of the
unit-sphere have been provided by Gilbert \cite{gilbert} in dimension three,
Hall \cite{hall} in any dimension when the unit-sphere is replaced by the
unit-cube and more recently by B\"urgisser, Cucker \& Lotz \cite{burgisser}. A
very minor consequence of Theorem 1.1 of this last work is the following: let 
$\widetilde{P}(f,n)$ denote  the probability to cover the unit-sphere with $n$ random 
circular caps which are independent, with uniformly-distributed centres and with a fractional area of $f$. Then
\begin{equation}\label{eq:500}
\lim_{n\to+\infty}\frac{1}{n}\log(1-\widetilde{P}(f,n))=\log(1-f).
\end{equation}
In particular, the fractional area occupied by a circular cap of angular
radius ${R}/{r}$ is
$$f_r=\frac{(d-1)\omega_{d-1}}{d\omega_{d}}\int_0^{R/r}\sin^{d-2}(\theta)\,\d\theta
\underset{r\to +\infty}{\sim}\frac{\omega_{d-1}R^{d-1}}{d\omega_{d} r^{d-1}}.$$
Consequently, we have
$$\P{{\mathfrak V}\geq r}\leq
e^{-\omega_d(r^d-R^d)}\sum_{n=0}^{+\infty}\frac{(\omega_d(r^d-R^d))^n}{n!}(1-\widetilde{P}(f_r,n))$$
and a direct application of (\ref{eq:500}) provides that 
$$\log\P{{\mathfrak V}\geq r}\underset{r\to +\infty}{\lesssim}-\omega_d r^d
f_r\underset{r\to+\infty}{\sim}-\frac{\omega_{d-1}}{d}R^{d-1}r.$$
This completes the proof of Proposition \ref{highdim}.\end{preuvesn}
\begin{remark}\label{visibilityinfinity}
Proposition \ref{highdim} implies that the total visibility ${\mathfrak V}$ is finite
almost surely. Nevertheless, when the
intensity measure of the underlying Poisson point process is of the form (in
spherical coordinates)
$r^{\alpha-1}\mathrm{d}r\mathrm{d}\sigma_d(u)$,
$\alpha\in{\mathbb R}$, 
it can be shown in the same way that the visibility at infinity
exists with positive probability as soon as $\alpha<1$.
\end{remark}

\section{Small or distant obstacles: convergence towards the law of extreme values}

When the size of the covering objects 
becomes smaller and the number of objects grows at the same time accordingly,
S. Janson \cite{janson} showed in very general setting that a particular
scaling yields a convergence towards the Gumbel law. We shall use this kind of result in two
contexts below:
\begin{itemize}
\item the asymptotics of the visibility when the obstacles are small (or equivalently when the intensity of the
centres is small);
\item the study of the visibility when there exists a large region around the origin with no obstacle at all.
\end{itemize}

\subsection{Small obstacles}

In this subsection, the radius $R$ of the obstacles will be deterministic but no longer
constant. For sake of clarity, we will denote by
${\mathfrak V}_R$ the visibility when the obstacles are discs of radius $R>0$. We
aim at giving the asymptotic behaviour of the visibility when
the size of the obstacles goes to $0$. Let us define the quantity
\begin{equation}
  \label{xir}
\xi_R=\omega_{d-1}R^{d-1}{\mathfrak
  V}_R+d(d-1)\log(R)-2(d-1)\log(-\log(R))-K_d
\end{equation}
where $$K_d=\log\left(\frac{d^{2(d-1)}(d-1)^{3(d-1)-1}\Gamma\left(\frac{d}{2}-\frac{1}{2}\right)^{2d-2}}{(d-1)!\pi^{\frac{(d-1)^2+1}{2}}2^{2d-3}\Gamma\left(\frac{d}{2}\right)^{d-2}}\right).$$
\begin{theorem}\label{smallR}
When $R$ goes to $0$, the quantity $\xi_R$ (provided by (\ref{xir})) converges in distribution to the extreme value
distribution, {\it i.e.\/} for every $u\in{\mathbb R}$, 
$$\lim_{R\to 0}\P{\xi_R\leq u}=\exp\left(-e^{-u}\right).$$
\end{theorem}
\begin{preuvesn}
The proof relies essentially on the application of a result due to Janson
(Lemma 8.1. in \cite{janson})
about random coverings of a compact Riemannian manifold by small geodesic
balls.

As before, we exploit the connection between the cumulative distribution
function of ${\mathfrak V}_R$ and the probability of covering the sphere with
circular caps:
\begin{eqnarray}
  \P{\xi_R\leq u}&=&{\mathbb P}\left( {\mathfrak V}_R\leq
    -\frac{d(d-1)}{\omega_{d-1}}\frac{\log(R)}{R^{d-1}}+\frac{2(d-1)}{\omega_{d-1}}
    \frac{\log(-\log(R))}{R^{d-1}}+\right.\nonumber\\
&&  \hspace*{1cm}  \left.\frac{K_d+u}{\omega_{d-1}}\frac{1}{R^{d-1}}\right)\nonumber\\
&=&{\mathbb P}(\mbox{the sphere of radius $f(R)$ is covered}\nonumber\\&&\hspace*{1cm} \mbox{by circular caps coming
    from the obstacles})\label{100}
\end{eqnarray}
where 
$$f(R)=-\frac{d(d-1)}{\omega_{d-1}}\frac{\log(R)}{R^{d-1}}+\frac{2(d-1)}{\omega_{d-1}}
\frac{\log(-\log(R))}{R^{d-1}}+\frac{K_d+u}{\omega_{d-1}}\frac{1}{R^{d-1}}.$$
Let us focus on this covering probability: the concerned obstacles are those
such that their centres $x$ are at distance $\rho\in (R,R+f(R))$. Their number is
Poisson distributed, of mean $\omega_d((R+f(R))^d-R^d)$. The set of all
$x/\|x\|$, where $x\in \xxx\cap [B(O,R+f(R))\setminus B(O,R)]$, is a
homogeneous Poisson
point process on the unit-sphere of intensity $({(f(R)+R)^d-R^d})/{d}$.

The induced
shadow of each of theses obstacles is a geodesic ball on the unit-sphere of
angular radius equal to $\arcsin(R/\rho)$ if $\rho\in (R,\sqrt{R^2+f(R)^2})$
and equal to $\arccos((f(R)^2+\rho^2-R^2)/2f(R)\rho)$ for $\rho\in
[\sqrt{R^2+f(R)^2},R+r)$ (see (\ref{length2})). In particular, it can be
verified that
the normalized geodesic radius $\Theta_R$ of this circular cap satisfies that 
\begin{equation}\label{convergencepsi}
\frac{1}{a_R}\Theta_R\overset{D}{\to}{\bf
  1}_{[1,+\infty)}(u)\frac{d}{u^{d+1}}\,\d u.
\end{equation}
where $a_R=\frac{R}{f(R)}\underset{R\to 0}{\to}0$.
Consequently, the required covering probability in (\ref{100}) is the probability that the
unit-sphere is covered by a Boolean model on the sphere of intensity
$\lambda_R=\frac{(f(R)+R)^d-R^d}{d}$, such that the geodesic balls have
i.i.d. radii distributed as $\Theta_R$ (with $\Theta_R$ satisfying the convergence
(\ref{convergencepsi})).

It only remains to verify that all the hypotheses of Lemma 8.1. in \cite{janson}
are satisfied (${\mathbb S}^{d-1}$ being a
  $(d-1)$-dimensional Riemannian manifold): 

\begin{itemize}
\item the
only notable difference is that we should not have $\frac{1}{a_R}\Theta_R$
converging in distribution but have it of fixed distribution for any
$R>0$. Nevertheless, the proof of Lemma 8.1. in \cite{janson} relies
essentially on convergence results [(7.15), (7.20), {\it ib.\/}] which also work in
this context without any changes;
\item the moments of order $((d-1)+\varepsilon)$ of the limit distribution
  obtained in (\ref{convergencepsi}) are finite for every $\varepsilon \in
  (0,1)$. Moreover, the moment of order $(d-1)$ is $d$ and the moment of order
  $(d-2)$ is $d/2$;
\item the constants $b$ and $\alpha$ defined in (\cite{janson}, Lemma 8.1) can be
  calculated:
$$b=\frac{\omega_{d-1}}{d\omega_d}\int_1^{+\infty}\frac{d}{u^{2}}\,\d u=\frac{\omega_{d-1}}{\omega_d}$$
and
$$\alpha=\frac{1}{d!}\left(\frac{\sqrt{\pi}\Gamma\left(\frac{d}{2}+1\right)}{\Gamma\left(\frac{d+1}{2}\right)}\right)^{d-1}\cdot\frac{\left(\frac{d}{2}\right)^{d-1}}{d^{d-2}};$$
\item The convergence (8.1) in \cite{janson} is satisfied:
\begin{eqnarray*}
&&\lim_{R\to   0}\left\{ba_R^{d-1}d\omega_d\lambda_R+\log(ba_R^{d-1})+\right.\\
&&\hspace*{1cm}\left.(d-1)\log(-\log(ba_R^{d-1}))-\log(\alpha)\right\}=u.
\end{eqnarray*}
Consequently, the proof of Theorem \ref{smallR} is complete.
\end{itemize}
\end{preuvesn}
\begin{remark} The same type of method and result should also occur in dimension two when the discs
are replaced by rotation-invariant i.i.d convex bodies.
\end{remark}
\begin{remark}
In any dimension, the result could be extended to radii of the form
$R_x=\varepsilon U_x$, $x\in \xxx$, where $\varepsilon$ goes to $0$ and the
$U_x$ are i.i.d.  bounded  random variables.
\end{remark}

\subsection{Conditioning by a large clearing}

We define here $S$ the clearing radius as 
$$S=\sup\{r>0;B_2(O,r)\subset {\mathbb R}^2\setminus{\mathfrak O}\}.$$
The distribution of $S$ is called the spherical
contact distribution. This section aims at estimating the distribution of the
maximal visibility ${\mathfrak V}$ conditionally on $S$. In particular we show that
when $S$ is large, ${\mathfrak V}$ is asymptotically
equivalent to $S$ (see Theorem \ref{clairiere}) and we estimate precisely the difference $({\mathfrak V}-S)$
via an extreme value result (see Theorem \ref{ConvMin}).

A first estimation based on techniques similar to the proofs of (\ref{eq:10}) and
(\ref{eq:10.1}) provides the following result:
\begin{theorem}\label{clairiere}
For every $\alpha\in (0,1)$, we have
$$\P{{\mathfrak V}\geq r+r^{\alpha}|S=r}\leq \P{{\mathfrak V}\geq r+r^{\alpha}|S\geq
r}=O\left(e^{-2E(R)r^{\alpha}}\right).$$
\end{theorem}
\begin{preuvesn}
Let us fix $r>0$. Conditionally on $\{S\geq r\}$, the process of couples
$(x,R_x)$ is a Poisson point process on ${\mathbb R}^2\times {\mathbb R}_+$ of
intensity measure ${\bf 1}_{\|x\|-R> r}dx\otimes \mu$. 

As in Proposition \ref{equiv-recouvre}, ${\mathfrak V}$ is greater than $r+u$, $u>0$, if
and only if the circle $C(O,r+u)$ is not totally hidden by the 'shadows' of
the obstacles. Moreover, the discs $B_2(x,R_x)$ which produce a non-empty
shadow are those which satisfy
$r< \|x\|-R_x < (r+u)$. The formula for the length of the shadow depends on
whether $\|x\|\leq \sqrt{(r+u)^2+R_x^2}$ or not (see equalities
(\ref{length2})). Consequently, if we only consider the shadows produced by the
discs $B_2(x,R_x)$ such that 
$$ (r+R_x<\|x\|<\sqrt{(r+u)^2+R_x^2})\;\;\mbox{ and }\;\;
u>\sqrt{r^2+2rR_x}-r,$$
then the probability of not covering the circle is greater. \\
When $u>\sqrt{r^2+2rR^*}-r$, the number of such discs is Poissonian, of mean
$\pi(u^2+2ru-2rE(R))$. Moreover, for these
discs, the length of the shadow decreases with $\|x\|$ and is minimal when
$\|x\|=\sqrt{(r+u)^2+R_x^2}$, equal to $L_{\mbox{\tiny{min}}}=({1}/{\pi})\arcsin
({R_x}/({(r+u)^2+R_x^2}))$ (see (\ref{length2})). In the sequel,
we denote by $\nu_{r,u}$ the distribution of $L_{\mbox{\tiny{min}}}$ and
$m_{r,u}$ its mean.

In conclusion, we have proved the following inequality: for every $u>\sqrt{r^2+2rR^*}-r$,
\begin{eqnarray}
&&\P{{\mathfrak V}\geq r+u|L\geq r}\nonumber\\
&&\leq e^{-\pi(u^2+2ru-2rE(R))}\sum_{n=0}^{+\infty}\frac{(\pi(u^2+2ru-2rE(R)))^n}{n!}(1-P(\nu_{r,u},n))\nonumber\\
&&\leq e^{-\pi(u^2+2ru-2rE(R))}\sum_{n=0}^{+\infty}\frac{(\pi(u^2+2ru-2rE(R)))^n}{n!}(1-P(\delta_{m_{r,u}},n)).\label{eq:100}
\end{eqnarray}
In particular, when $u=r^{\alpha}$, $0<\alpha<1$, we have
$$m_{r,r^{\alpha}}\underset{r\to +\infty}{\sim}\frac{1}{\pi}\frac{E(R)}{r}\;\;\mbox{
and } \pi(u^2+2ru-2rE(R))\underset{r\to +\infty}{\sim}2\pi r^{1+\alpha}$$
Using the inequality (\ref{eq:sheppameliore}) and inserting the two previous estimates in (\ref{eq:100}), we obtain
the required result, i. e. 
$$\P{{\mathfrak V}\geq r+r^{\alpha}|S\geq
r}=O\left(e^{-2E(R)r^{\alpha}}\right).$$
Finally, it remains to study the distribution of ${\mathfrak V}$ conditionally
on ${S=r}$. We remark that conditionally on $\{S=
r\}$, the process is the same as in the case of the conditioning on $\{S\geq
r\}$ with a supplementary random disc $B_2(Y,R_Y)$
such that $R_Y$ is $\mu$-distributed and conditionally on $R_Y$, $Y$ is
uniformly distributed on the circle $C(O,r+R_Y)$. Since there is one more
obstacle, the
maximal visibility must be lesser than in the case of the conditioning on
$\{S\geq r\}$.
\end{preuvesn}
Theorem \ref{clairiere} implies that the difference $({\mathfrak V}-S)$ is
negligible in front of $S$ but we can get a far more precise three-terms
development in the following way: for every $r>0$, we denote by ${\mathfrak V}_r$ a random variable distributed as ${\mathfrak V}$ when the Boolean model is conditioned on $\{S\ge r\}$, {\it i.e.\/} on not having any grain at distance lesser than $r$ from the origin. Let us define the quantity
\begin{equation}
  \label{psil}
\psi_r=\omega_{d-1}E_{\mu}(R^{d-1})({\mathfrak V}_r-r)-(d-1)\log(r)-(d-1)\log(\log(r))-K_d'
\end{equation}
where 
\begin{eqnarray*}
&&K_d'=\log\left(\frac{1}{(d-1)!}\left(\frac{\sqrt{\pi}\Gamma\left(\frac{d+1}{2}\right)}{\Gamma\left(\frac{d}{2}\right)}\right)^{d-2}\frac{E_{\mu}(R^{d-2})^{d-1}}{E_{\mu}(R^{d-1})^{d-2}}\right)\\
&&\hspace*{2cm}+(d-1)\log(d-1)-\log\left(\frac{\omega_{d-1}E_{\mu}(R^{d-1})}{d\omega_d}\right).
  \end{eqnarray*}
For every $t\in{\mathbb R}$, we have
\begin{eqnarray*}
\P{\psi_{r}\le t}&=&{\mathbb P}\left(\omega_{d-1}E_{\mu}(R^{d-1})({\mathfrak V}-r)-\right.\\
&&\hspace*{1cm}\left.(d-1)\log(r)-(d-1)\log(\log(r))-K\le t|S\ge r\right).  
\end{eqnarray*}
\begin{theorem}\label{ConvMin}
When $r$ goes to $\infty$, the quantity $\psi_r$ (provided by (\ref{psil})) converges in distribution to the extreme value
distribution, {\it i.e.\/} for every $t\in{\mathbb R}$, 
$$\lim_{r\to +\infty}\P{\psi_{r}\le t}=\exp\left(-e^{-t}\right).$$
\end{theorem}
\begin{preuvesn} 
As previously in the case of small obstacles, the proof relies essentially on the application of a result due to Janson
(Lemma 8.1. in \cite{janson})
about random coverings of a compact Riemannian manifold by small geodesic
balls. Let us consider the quantity
$$f(r)=\frac{d-1}{\omega_{d-1}E_{\mu}(R^{d-1})}\log(r)+
\frac{d-1}{\omega_{d-1}E_{\mu}(R^{d-1})}\log(\log(r))+\frac{K+t}{\omega_{d-1}E_{\mu}(R^{d-1})}.$$  
such that $(\psi_r\le t)\Longleftrightarrow ({\mathfrak V}_r-r\le
f(r))$. The connection with a covering probability is the following:
\begin{eqnarray*}
\P{\psi_r\le u}&=&{\mathbb P}(\mbox{the sphere of radius $(r+f(r))$ is covered}\\
&&\mbox{~~~~~~~~by circular caps coming from the obstacles}).
\end{eqnarray*}
It remains to investigate asymptotics of this covering probability: the concerned obstacles are those
such that their centers $x$ are at distance $\rho\in (r+R,r+f(r)+R)$. Their number is
Poisson distributed, of mean $\omega_dE_{\mu}[(R+f(r)+r)^d-(R+r)^d]$. 

The induced
shadow of each of theses obstacles is a geodesic ball on the unit-sphere of
angular radius equal to:
\begin{itemize}
\item $\arcsin(R/\rho)$ if $\rho\in (R+r,\sqrt{R^2+(f(r)+r)^2})$,
\item $\arccos(((f(r)+r)^2+\rho^2-R^2)/2(f(r)+r)\rho)$ if $\rho\geq \sqrt{R^2+(f(r)+r)^2}$ and
$\rho<R+r+f(r))$ (see (\ref{length2})). 
\end{itemize}
In particular, it can be verified that
the normalized geodesic radius $\Theta_r$ of this circular cap satisfies that 
$r\Theta_r\overset{D}{\to}\mu.$ In the rest of the proof, we will use the quantity 
$a_r={1}/{r}$ in order to be as close as possible to the notations of Janson's lemma.

Consequently, the required covering probability is the probability that the
unit-sphere is covered by a Boolean model on the sphere of intensity
$$\lambda_r=\frac{1}{d}E_{\mu}[(R+f(r)+r)^d-(R+r)^d]\underset{r\to +\infty}{\sim} r^{d-1}f(r),$$ such that the geodesic balls have
i.i.d. radii distributed as $\Theta_r$.

As in the proof of Theorem \ref{smallR}, we verify that all the hypotheses of Lemma 8.1. in \cite{janson}
are satisfied:
\begin{itemize}
\item the moments of order $((d-1)+\varepsilon)$ of the limit distribution $\mu$ are finite for every $\epsilon>0$. 
\item the constants $b$ and $\alpha$ defined in (\cite{janson}, Lemma 8.1) can be
  calculated:
$$b=\frac{\omega_{d-1}E_{\mu}(R^{d-1})}{d\omega_d}$$
and
$$\alpha=\frac{1}{(d-1)!}\left(\frac{\sqrt{\pi}\Gamma\left(\frac{d+1}{2}\right)}{\Gamma\left(\frac{d}{2}\right)}\right)^{d-2}\cdot\frac{E_{\mu}(R^{d-2})^{d-1}}{E_{\mu}(R^{d-1})^{d-2}};$$
\item The convergence (8.1) in \cite{janson} is satisfied:
\begin{eqnarray*}
&&\lim_{r\to  +\infty}\left\{ba_r^{d-1}d\omega_d\lambda_r+\log(ba_r^{d-1})+\right.\\
&&\hspace*{1cm} \left.(d-1)\log(-\log(ba_r^{d-1}))-\log(\alpha)\right\}=t.
\end{eqnarray*}
\end{itemize}
Consequently, the proof of Theorem \ref{smallR} is complete.
\end{preuvesn}
\section*{Appendix 1: proof of the lower bound via direct computation, estimate (\ref{mino-main})}
It is quite reasonable to try to obtain directly a lower bound on the tail probability, let us explain the sketch of the proof:
the visibility is greater than $r$ if and only if there exists a direction in which one can see farther, so that if one 
discretises the circle $\partial B_{2}(0,r)$, one could argue that there exists one of those directions such that the visibility
is greater than $r$, the number of those directions is $2/\pi r$, hence the order $r\,\exp(-2rR)$. We shall make this 
statement more rigorous below.\\

Let us take $\zeta\in(0, 2/R)$, we define $N_r$ as the integer part of $\zeta s$ and $\theta_r={2\pi}/({\zeta r}),$
and define the points $A_{k,r} = (r,k\theta_r)$ for $k\in\{0,\ldots,N_r-1\}$. We see easily that if we
define
$$G_{k,r}=(B_{2}(0,R)\oplus [0,A_{k,r}])\setminus B_{2}(0,R),$$
then for $r>R$ and $k_{1}\neq k_{2}$ one has
$$G_{k_{1},r}\cap G_{k_{2},r}\subset B_{2}(0,r).$$
The sets $G_{k,r}$ are called `fingers', we denote by 
$\rho_r=R/\sin(\theta_{r}/2)$ the maximal norm of a point belonging to the 
intersection of two fingers, we shall denote by $E_{k,r}$ the intersection $G_{k,r}\cap B_{2}(0,\rho_{r})$
and $F_{k,r}=G_{k,r}\setminus E_{k,r}$.\\

\begin{figure}[h]
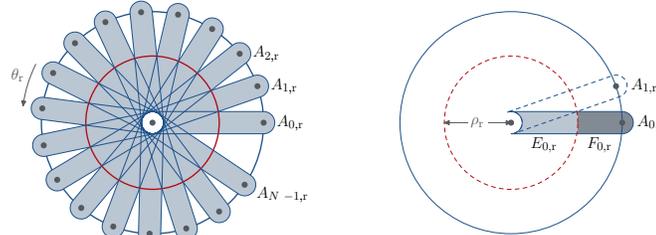

\begin{center}
{\white \begin{tabular}{c}
1\\
1\\
1\\
1\\
1\\
1\\
\end{tabular}}
{\white rien$\!\!\!\!\!\!\!\!\!$}\includegraphics{minoDebutFleur.0}
{\white ~~ ~~~~~~~~~~~~~~~~~~~~ ~~~ ~~~ ~~~~~ ~~~~}
\includegraphics{minoDefEF.0}
\end{center}
\caption{The points $A_{i,r}$ and their associated `fingers'.
The point $A_{k,r}$ is visible if and only if no point of the
process $\xxx$ belongs to $G_{k,r}=E_{k,r}\cup F_{k,r}$.}\label{roue-dentee}
\end{figure}

We will assume from now on that $r$ is large enough.
If at least one of those points $A_{k,r}$ is not shadowed by the
discs intersecting $B_{2}(0,r)$, the visibility ${\mathfrak V}$ is greater than $r$: hence the probability
of this event is greater than the probability that one of the `fingers' $G_{k,r}$ in figure \ref{roue-dentee}
does not contain a point of $\xxx$.

\begin{figure}[h*]
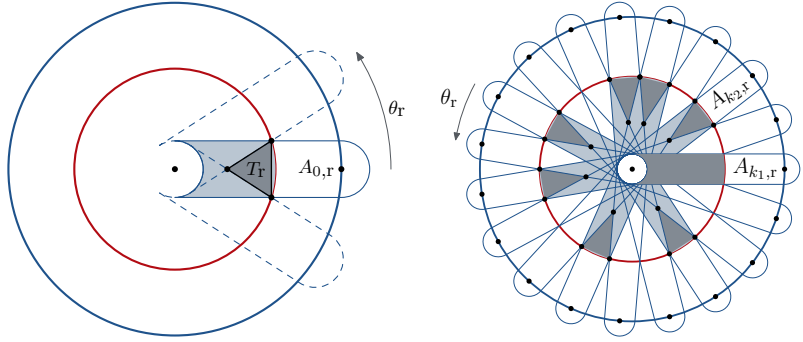

\begin{center}
\begin{tabular}{cc}
\includegraphics[height=4.5cm]{minoDefT0.0}&
\includegraphics[height=4.5cm]{minoUnionEki.0}\\
\end{tabular}
\caption{Left: definition of $T_r$. Right: 
the set $E_{k_1,r}\cup\cdots\cup E_{k_{\ell},r}$ in light grey and in dark grey
the subset with area $\aire{2}{E_{0,r}}+(\ell-1)\aire{2}{T_r}$.}
\label{fig:DefTzero}
\end{center}
\end{figure}

We have
$$\aire{2}{F_{k,r}} = 2R(1-\kappa)s+\O(1),$$
where $\kappa=R\zeta \pi^{-1}$, so that $1-\kappa>0$,  and $\O(1)$ means a bounded function. 
Let us define the following events
\begin{eqnarray*}
V_{k,r}&=&\left\{ {\xxx}\cap G_{k,r}=\vide \right\},\\
Z_{k,r}&=&\left\{ {\xxx}\cap  E_{k,r}=\vide \right\},\\
W_{k,r}&=&\left\{ {\xxx}\cap F_{k,r}=\vide \right\}.
\end{eqnarray*}
We have
$$V_{k,r}=Z_{k,r}\cap W_{k,r},$$
and we want to evaluate $\P{\bigcup_{k=0}^{N_r-1}V_{k,r}}$, this is 
equal thanks to Poincar\'e's formula to
\begin{equation}
  \label{poincare}
\P{\bigcup_{k=0}^{N_r-1}V_{k,r}}
= \sum_{k=0}^{N_r-1}\P{V_{k,r}}+\sum_{\ell=2}^{N_r}(-1)^{\ell-1}
\sum_{k_1<\cdots<k_\ell}\P{V_{k_1,r}\cap\cdots\cap V_{k_\ell,r}}.
  \end{equation}
We shall prove that the dominating term in this expansion is the first one: it rewrites as
\begin{eqnarray}
\sum_{k=0}^{N_r-1}\P{V_{k,r}} &=& \displaystyle\sum_{k=0}^{N_r-1}
\fexp{-\aire{2}{E_{k,r}\cup F_{k,r}}}\nonumber\\
&=& N_r \fexp{-2Rr}.\nonumber
\end{eqnarray}
Let us consider $\ell\geq 2$ and $0\leq k_1<\cdots <k_\ell<N_{r}$, we have:
$$
\P{V_{k_1,r}\cap\cdots\cap V_{k_\ell,r}} 
= \P{\bigcap_{i=1}^\ell Z_{k_i,r}\cap \bigcap_{i=1}^{\ell}W_{k_i,r}},
$$
where the events $\bigcap_{i=1}^\ell Z_{k_i,r}$, $W_{k_1,r}$, \ldots, $W_{k_\ell,r}$ 
are independent as the sets $\bigcup_{i=1}^\ell E_{k_i,r}$ and $F_{k_1,r},\ldots, F_{k_\ell,r}$ 
are disjoint, hence:
\begin{eqnarray}
\P{V_{k_1,r}\cap\cdots\cap V_{k_\ell,r}} 
&=&\fexp{-\aire{2}{\bigcup_{i=1}^\ell E_{k_i,r}}}\fexp{-\ell\aire{2}{F_{0,r}}}.\nonumber
\end{eqnarray}
\noindent To estimate this probability, let us introduce
the triangle $T_r$ which is the greatest triangle included in $E_{0,r}\backslash\left(\cup_{i=1}^{N_r-1}E_{k,r}\right)$
(see \ref{fig:DefTzero}). We easily get:
$$\aire{2}{T_r}= 2R\,\frac{\kappa}{4}s+\O(1).$$
For each $k$, $E_{k,r}$ contains a triangle that is isometric to $T_r$ and disjoint from all others $E_{k',r}$, hence we have:
\begin{eqnarray*}
\aire{2}{\bigcup_{i=1}^\ell E_{k_i,r}}
&\geq&(\ell-1)\aire{2}{T_r}+\aire{2}{E_{0,r}},\\
&\geq&\ell\aire{2}{T_r}+\aire{2}{E_{0,r}}-\aire{2}{T_{r}}.
\end{eqnarray*}

Thus for all choice of $\ell\geq 2$ and $0\leq k_1<k_2<\cdots<k_{\ell}<N_{r}$ we have:
\begin{eqnarray*}
&&\fexp{-\aire{2}{B_2(0,\rho_r)}}\fexp{-\ell\aire{2}{F_{0,r}}}\leq \P{\cap_{i=1}^\ell V_{k_i,r}}\leq\\
&&\fexp{-\aire{2}{E_{0,r}}+\aire{2}{T_{r}}}\fexp{-\ell(\aire{2}{F_{0,r}}+\aire{2}{T_r})}.
\end{eqnarray*}
The number of such terms is $\binom{N_r}{\ell}$, hence the sum 
$S_{\ell}$ of all those terms
satisfies
$$\left|S_{\ell}\right|\leq\fexp{\aire{2}{T_r}-\aire{2}{E_{0,r}}}
\binom{N_r}{\ell}\fexp{-\ell(\aire{2}{F_{0,r}}+\aire{2}{T_r})}$$
and the residual term $S=\sum_{\ell=2}^{N_r}(-1)^{\ell-1}S_\ell$ is bounded from above by:
\begin{eqnarray}
|S|&\leq&\fexp{\aire{2}{T_r}-\aire{2}{E_{0,r}}}\times\nonumber\\
&&\sum_{\ell=2}^{N_r}\binom{N_r}{\ell}\fexp{-\ell(\aire{2}{F_{0,r}}+\aire{2}{T_r})},\nonumber\\
&\leq&\fexp{\aire{2}{T_r}-\aire{2}{E_{0,r}}}\times\nonumber\\
&&\Bigl(
\bigl(1+\fexp{-(\aire{2}{F_{0,r}}+\aire{2}{T_r})}\bigr)^{N_r}\nonumber\\
&&-1-N_r\fexp{-(\aire{2}{F_{0,r}}+\aire{2}{T_r})}
\Bigr),\nonumber\\
&\leq&\frac{N_r^2}{2}\, \fexp{\aire{2}{T_r}-\aire{2}{E_{0,r}}}\times\nonumber\\
&&\fexp{-2(\aire{2}{F_{0,r}}+\aire{2}{T_r})}\left(1+\O(1)\right).
\end{eqnarray}
\def\e{\textrm{e}}
Using the asymptotic expansions of $F_{0,r}$, $E_{0,r}$ and $T_{r}$ and the unpper bound $N_r\leq \zeta r$, we obtain:
\begin{eqnarray*}
|S|&\leq& \frac{1}{2}\,
\fexp{-\frac{3}{2}\,R\kappa s+\O(1)}\zeta^2 r^2\times\\
&&\fexp{-2R\left(2-\frac{3}{2}\kappa\right)r+\O(1)}\left(1+\O(1)\right)\\
&\leq&C
\zeta^2 r^2\fexp{-2Rr\left(1+\left(1-\frac{3}{4}\kappa\right)\right)}.
\end{eqnarray*}
hence $S=\o(\zeta r\exp(-2Rr))$, which completes the proof of estimate \ref{mino-main}.

\begin{remark}
In the proof above, we could have taken only one triangle to obtain that the sum $S$
is negligible with respect to the first term $\zeta s\fexp{-2Rr}$, however the
accuracy of the development would have been less interesting. Discs with bounded random radius $R\in[\Rdown,\Rup]$
can also be treated this way, at a cost of a loss on the accuracy because of a non-optimal size of the fingers.
\end{remark}
\section*{Appendix 2: proof of Lemma \ref{geometriclemma}}
For sake of simplicity, we call $x=r{\bf u}$. Let us denote by $y$ (resp. $y'$) a point in the intersection of $(r{\bf
  u}+L)$  with its tangent line emanating from $O$ and situated on the
left-hand side (resp. on the right-hand side) of the half-line $(O+{\mathbb
  R}_+{\bf u})$. We define $z$ (resp. $z'$) as the orthogonal projection of
$y$ (resp. $y'$) on $(O+{\mathbb R}_+{\bf u})$ and $\alpha$ (resp. $\alpha'$)
as the angle between $(O+{\mathbb R}_+{\bf u})$ and $(O+{\mathbb R}_+ y)$
(resp. $(O+{\mathbb R}_+ y')$). Then we have
\begin{equation}\label{angle}
\Psi(r{\bf
  u}+L)=\alpha+\alpha'=\arctan\left(\frac{||y-z||}{||z||}\right)+\arctan\left(\frac{||y'-z'||}{||z'||}\right).
\end{equation}
Let us now describe $W_u(L)$: there exist two points $w$ and $w'$ ($w$ being
on the left-hand side of $(O+{\mathbb R}_+{\bf u})$) such that
\begin{equation}\label{diam}
W_u(L)=\mbox{dist}(w,O+{\mathbb R}{\bf u})+\mbox{dist}(w',O+{\mathbb R}{\bf
u})
\end{equation}
where $\mbox{dist}(\cdot,O+{\mathbb R}{\bf u})$ is the Euclidean distance to
the line $(O+{\mathbb R}{\bf u})$.
Comparing (\ref{angle}) with (\ref{diam}), we observe that we only have to prove that 
\begin{equation}
  \label{eq:105}
\lim_{r\to +\infty} r\alpha=\lim_{r\to +\infty}
r\arctan\left(\frac{||y-z||}{||z||}\right)=\mbox{dist}(w,O+{\mathbb R}{\bf
  u})
\end{equation}
and 
\begin{equation}
\lim_{r\to +\infty}r\alpha'=\lim_{r\to +\infty}
r\arctan\left(\frac{||y'-z'||}{||z'||}\right)=\mbox{dist}(w',O+{\mathbb R}{\bf
  u}),
  \end{equation}
the limits being uniform as required. Let us now concentrate on the first
limit (the second can be proved in the same way):

since $||y-z||\le ||y-x||\le D$ and $||z-x||\le ||y-x||\le D$, we have for
every $r\ge D$
\begin{equation}
  \label{eq:tangente}
\tan(\alpha)=\frac{||y-z||}{||z||}\le\frac{||y-z||}{r-D}\le \frac{D}{r-D}.
  \end{equation}
Moreover, 
$$0\le \mbox{dist}(w,O+{\mathbb R}{\bf u})-||y-z||=||y-w||\sin(\beta)\le D\sin(\beta)$$
where $\beta$ is the angle between $(O+{\mathbb R}_+{\bf u})$ and the line
from $y$ to $w$. Since $y$ is a contact point of a support line of $L$ and $w$ is in $L$, this angle $\beta$ must
necessarily be lesser than $\alpha$. Consequently, we get by a direct use of
(\ref{eq:tangente}) that
\begin{equation}
  \label{eq:104}
  0\le \mbox{dist}(w,O+{\mathbb R}{\bf u})-||y-z||\le D\sin(\alpha)\le
  D\tan(\alpha)\le \frac{D^2}{r-D}.
\end{equation}
Inserting this last estimate in the first equality of (\ref{eq:tangente}), we have
$$r\arctan\left(\frac{\mbox{dist}(w,O+{\mathbb R}{\bf
      u})-\frac{D^2}{r-D}}{r+D}\right) \le r\alpha\le
r\arctan\left(\frac{\mbox{dist}(w,O+{\mathbb R}{\bf u})}{r-D}\right),$$
which provides the required convergence result (\ref{eq:105}) with a
uniformity in ${\bf u}$ and in $L$.



\begin{thebibliography}{10}

\bibitem{ballani}
{\sc Ballani, F.} (2006).
\newblock On second-order characteristics of germ-grain models with convex
  grains.
\newblock {\em Mathematika\/} {\bf 53,} 255--285 (2007).

\bibitem{schramm}
{\sc Benjamini, I., Jonasson, J., Schramm, O. and Tykesson, J.}
\newblock Visibility to infinity in the hyperbolic plane, despite obstacles.
\newblock http://arxiv.org/abs/0807.3308 2008.

\bibitem{burgisser}
{\sc B\"urgisser, P., Cucker, F. and Lotz, M.}
\newblock Coverage processes on spheres and condition numbers for linear
  programming.
\newblock http://arxiv.org/abs/0712.2816 2008.

\bibitem{calka03}
{\sc Calka, P.} (2002).
\newblock The distributions of the smallest discs containing the
  {P}oisson-{V}oronoi typical cell and the {C}rofton cell in the plane.
\newblock {\em Adv. in Appl. Probab.\/} {\bf 34,} 702--717.

\bibitem{CMPB09}
{\sc Calka, P., Michel, J. and Porret-Blanc, S.} (2009).
\newblock Visibilit\'e dans le mod\`ele bool\'een.
\newblock To appear in {\em C. R. Math. Acad. Sci. Paris\/}.

\bibitem{capassovilla}
{\sc Capasso, V. and Villa, E.} (2005).
\newblock Survival functions and contact distribution functions for
  inhomogeneous, stochastic geometric marked point processes.
\newblock {\em Stoch. Anal. Appl.\/} {\bf 23,} 79--96.

\bibitem{dvoretzky}
{\sc Dvoretzky, A.} (1956).
\newblock On covering a circle by randomly placed arcs.
\newblock {\em Proc. Nat. Acad. Sci. U.S.A.\/} {\bf 42,} 199--203.

\bibitem{fan}
{\sc Fan, A.-H. and Wu, J.} (2004).
\newblock On the covering by small random intervals.
\newblock {\em Ann. Inst. H. Poincar\'e Probab. Statist.\/} {\bf 40,} 125--131.

\bibitem{gilbert}
{\sc Gilbert, E.~N.} (1965).
\newblock The probability of covering a sphere with {$N$} circular caps.
\newblock {\em Biometrika\/} {\bf 52,} 323--330.

\bibitem{hall}
{\sc Hall, P.} (1985).
\newblock On the coverage of {$k$}-dimensional space by {$k$}-dimensional
  spheres.
\newblock {\em Ann. Probab.\/} {\bf 13(3),} 991--1002.

\bibitem{heinrich}
{\sc Heinrich, L.} (1998).
\newblock Contact and chord length distribution of a stationary {V}orono\u\i\
  tessellation.
\newblock {\em Adv. in Appl. Probab.\/} {\bf 30,} 603--618.

\bibitem{jankovic}
{\sc Jankovi{\'c}, V.} (1996).
\newblock Solution of one problem of {G}. {P}\'olya.
\newblock {\em Mat. Vesnik\/} {\bf 48,} 47--50.

\bibitem{janson}
{\sc Janson, S.} (1986).
\newblock Random coverings in several dimensions.
\newblock {\em Acta Math.\/} {\bf 156,} 83--118.

\bibitem{kahane}
{\sc Kahane, J.-P.} (1990).
\newblock Recouvrements al\'eatoires et th\'eorie du potentiel.
\newblock {\em Colloq. Math.\/} {\bf 60/61,} 387--411.

\bibitem{lastschass}
{\sc Last, G. and Schassberger, R.} (2001).
\newblock On the second derivative of the spherical contact distribution
  function of smooth grain models.
\newblock {\em Probab. Theory Related Fields\/} {\bf 121,} 49--72.

\bibitem{molchanov}
{\sc Molchanov, I.} (2005).
\newblock {\em Theory of random sets}.
\newblock Probability and its Applications (New York). Springer-Verlag London
  Ltd., London.

\bibitem{muchestoyan}
{\sc Muche, L. and Stoyan, D.} (1992).
\newblock Contact and chord length distributions of the {P}oisson {V}orono\u\i\
  tessellation.
\newblock {\em J. Appl. Probab.\/} {\bf 29,} 467--471.

\bibitem{mullerstoyan}
{\sc M{\"u}ller, A. and Stoyan, D.} (2002).
\newblock {\em Comparison methods for stochastic models and risks}.
\newblock Wiley Series in Probability and Statistics. John Wiley \& Sons Ltd.,
  Chichester.

\bibitem{polya}
{\sc P\'olya, G.} (1918).
\newblock Zahlentheoretisches und {W}ahrscheinlichkeitstheoretisches {\"u}ber
  die {S}ichweite im {W}alde.
\newblock {\em Archiv der {M}athematik und {P}hysik, ser. 3\/} 135--142.

\bibitem{ratajsaxl}
{\sc Rataj, J. and Saxl, I.} (1997).
\newblock Boolean cluster models: mean cluster dilations and spherical contact
  distances.
\newblock {\em Math. Bohem.\/} {\bf 122,} 21--36.

\bibitem{shepp}
{\sc Shepp, L.~A.} (1972).
\newblock Covering the circle with random arcs.
\newblock {\em Israel J. Math.\/} {\bf 11,} 328--345.

\bibitem{siegelholst}
{\sc Siegel, A.~F. and Holst, L.} (1982).
\newblock Covering the circle with random arcs of random sizes.
\newblock {\em J. Appl. Probab.\/} {\bf 19,} 373--381.

\bibitem{stevens}
{\sc Stevens, W.~L.} (1939).
\newblock Solution to a geometrical problem in probability.
\newblock {\em Ann. Eugenics\/} {\bf 9,} 315--320.

\bibitem{skm}
{\sc Stoyan, D., Kendall, W.~S. and Mecke, J.} (1987).
\newblock {\em Stochastic geometry and its applications}.
\newblock Wiley Series in Probability and Mathematical Statistics: Applied
  Probability and Statistics. John Wiley \& Sons Ltd., Chichester.
\newblock With a foreword by D. G. Kendall.

\bibitem{yadin}
{\sc Yadin, M. and Zacks, S.} (1982).
\newblock Random coverage of a circle with applications to a shadowing problem.
\newblock {\em J. Appl. Probab.\/} {\bf 19,} 562--577.

\bibitem{zacksLN}
{\sc Zacks, S.} (1994).
\newblock {\em Stochastic visibility in random fields} vol.~95 of {\em Lecture
  Notes in Statistics}.
\newblock Springer-Verlag, New York.
\newblock With 1 IBM-PC floppy disc (3.5 inch; HD).

\end{thebibliography}
\end{document}